\documentclass[11pt]{article}

\usepackage{titlesec}

\titleformat{\paragraph}
{\normalfont\normalsize\bfseries}{\theparagraph}{1em}{}
\titlespacing*{\paragraph}
{0pt}{3.25ex plus 1ex minus .2ex}{1.5ex plus .2ex}

\usepackage{mathrsfs}

\usepackage[dvips]{epsfig}
\usepackage{graphicx}
\usepackage{amssymb,graphics,amsmath,amsthm,amsopn,amstext,amsfonts, algorithm,algorithmic}
\usepackage{color}

\newcommand{\wt}{\widetilde}

\usepackage[left=3cm,right=2cm,top=3cm,bottom=2cm]{geometry}

\newcommand{\tr}{ \mbox{\textbf{Tr}} }

\newcommand{\Argmin}{\mbox{Argmin}}

\newcommand{\Ical}{\mathcal I}

\newcommand{\Lcal}{\mathcal L}

\newcommand{\rank}{\mbox{rank}}

\newcommand{\mycut}[1]{{}}

\newtheorem{theorem}{Theorem}[section] 
\newtheorem{lemma}{Lemma}[section] 

\usepackage[utf8]{inputenc}
\usepackage[T1]{fontenc}
\usepackage{authblk}
\makeatletter
\def\thanks#1{\protected@xdef\@thanks{\@thanks
        \protect\footnotetext{#1}}}
\makeatother

\usepackage{blindtext}
\usepackage{multicol}

\begin{document}
\title{\textbf{\color{black}{Statistical Proxy based Mean-Reverting Portfolios with Sparsity and Volatility Constraints}
}}
\author[1]{Ahmad Mousavi$^*$\thanks{Emails: mousavi@american.edu and  gmichail@ufl.edu
}} 
\author[2]{George Michailidis\thanks{}}
\affil[1]{Department of Mathematics and Statistics, American University}
\affil[2]{Department of Statistics, University of Florida}
\date{}
\maketitle
\begin{abstract}
Mean-reverting portfolios with volatility and sparsity constraints are of prime interest to practitioners in finance since they are both profitable and well-diversified, while also managing risk and minimizing transaction costs. Three main measures that serve as statistical proxies to capture the mean-reversion property are predictability, portmanteau criterion, and crossing statistics. If in addition, reasonable volatility and sparsity for the portfolio are desired, a convex quadratic or quartic objective function, subject to nonconvex quadratic and cardinality constraints needs to be minimized. In this paper, we introduce and investigate a comprehensive modeling framework that incorporates all the previous proxies proposed in the literature and develop an effective \textit{unifying}  algorithm that is enabled to obtain a Karush–Kuhn–Tucker (KKT) point under mild regularity conditions. Specifically, we present a tailored penalty decomposition method that approximately solves a sequence of penalized subproblems by a block coordinate descent algorithm. To the best of our knowledge, our proposed algorithm is the first method for directly solving volatile, sparse, and mean-reverting portfolio problems based on the portmanteau criterion and crossing statistics proxies. 
Further, we establish that the convergence analysis can be extended to a nonconvex objective function case if the starting penalty parameter is larger than a finite bound and the objective function has a bounded level set. Numerical experiments on the S\&P 500 data set, demonstrate the efficiency of the proposed algorithm in comparison to a semidefinite relaxation-based approach and suggest that the crossing statistics proxy yields more desirable portfolios. 
\end{abstract}


\section{Introduction} 

Portfolios exhibiting mean-reverting behavior are of great interest to practitioners due to their predictability property that creates arbitrage opportunities. The task is to construct a portfolio that tends to return to its average value over time and use its performance to make trading decisions. 
Traditionally such portfolios are obtained by considering a linear combination of assets that are stationary through co-integration methods. However, such an approach often results in a large number of assets, which may not own reasonable volatility; making such portfolios practically inapplicable \cite{
fogarasi2013sparse,zhang2020sparse}. Therefore, along with mean-reversion, reasonable volatility is a mandatory property in realistic situations. Moreover, trading costs can be notably reduced by limiting the number of assets, and hence sparse portfolios offer advantages (note that the notion of sparsity is widely used in signal processing and machine learning applications; e.g., \cite{mousavi2020survey} and references therein). 
{\color{black}
Hence, constructing a mean-reverting portfolio exhibiting reasonable volatility and comprising a relatively smaller number of assets has found attention in the literature  \cite{mousavi2022penalty, cuturi2016mean}.} 
Recognize that the concept of sparsity holds significant relevance across various applications within contemporary science \cite{mousavi2023cardinality, shen2018least, moosaei2023sparse, ayanzadeh2019quantum}.

Three standard proxies that reflect the mean-reversion property are the predictability measure, the portmanteau criterion, and crossing statistics \cite{cuturi2013mean}.
The predictability measure, initially defined for stationary processes in \cite{box1977canonical} and later generalized for non-stationary processes by \cite{bewley1994comparison}, assesses how similar a time series is to a noise process.
{\color{black} Let $A_i$'s for $i=0,1,\dots, q$ be empirical autocovariance matrices related to a multivariate time series (see (\ref{def: autocov}) for details), then}
the mean-reversion property via the predictability proxy can be captured using the following problem \cite{cuturi2013mean}:
\begin{equation}\label{pr: pred}
\begin{aligned}
\min_{x\in \mathbb{R}^n}\quad 
x^TA_1x \qquad 
\textrm{subject to} \qquad 
x^Tx=1.
\end{aligned}
\end{equation}
The portmanteau statistic, proposed by \cite{ljung1978measure}, is used to determine if a process is a white noise. 
{\color{black}
In particular, mean reversion is obtained via portamento criterion proxy as follows \cite{cuturi2013mean}:
\begin{equation}\label{pr: port}
\begin{aligned}
\min_{x\in \mathbb{R}^n}\quad 
\sum_{i=2}^q \left(x^TA_ix\right)^2  
\qquad \textrm{subject to} \qquad 
x^Tx=1.
\end{aligned}
\end{equation}
}
The choice of $q$ depends on factors such as data availability, portfolio size, risk tolerance, and computational complexity; hence, in applications, sensitivity analysis is recommended for selecting an appropriate $q$  \cite{du2022high,imai2020statistical}.
The crossing statistic calculates the likelihood of a process crossing its mean per unit of time \cite{ylvisaker1965expected}. 
{\color{black}
Precisely, the crossing statistics proxy results in mean reversion by \cite{cuturi2013mean}:
\begin{equation}\label{pr: cross}
\begin{aligned}
\min_{x\in \mathbb{R}^n}\quad 
 x^TA_1x+\gamma \sum_{i=2}^q \left(x^TA_ix\right)^2  
\qquad \textrm{subject to} \qquad 
x^Tx=1,
\end{aligned}
\end{equation}
$\gamma$ serves as a  hyperparameter, governing the influence of lag-$i$ autocovariance matrices (see (\ref{def: autocov})) on the mean reversion characteristic of the target portfolio.
On the other hand, recall that sufficient volatility can be obtained by imposing a constraint on the variance, that is,  $xA_0x \ge \phi$ for some positive threshold $\phi$, where $A_0$ is the covariance matrix. Hence, statistical proxy-based problems (\ref{pr: pred}), (\ref{pr: port}), and (\ref{pr: cross}) together with volatility and sparsity constraints can be formulated into the following comprehensive model:
\begin{equation}\label{pr: p-original} \tag{\mbox{$P$}}
\begin{aligned}
\min_{x\in \mathbb{R}^n}\quad & \ 
 f(x):=\alpha x^TA_1x+\gamma \sum_{i=2}^q \left(x^TA_ix\right)^2  
\\ \textrm{subject to} \quad & 
 x^TA_0x\ge \phi, \quad 
x^Tx=1,    \quad 
\text{and} \quad    
\|x\|_0\le k,   
\end{aligned}
\end{equation}
where  $\alpha \in \{0,1\}, 2\le q\in \mathbb N, A_0\succ 0, A_i\succeq 0$   {\color{black} for $i=1,\dots, q, \gamma\ge 0,\phi >0$, and $k\le n$ (for further details see \ref{params})). We emphasize that} autocovariance matrices provide information about the temporal dependence structure among variables in a time series data set and are commonly used in time series analysis and econometrics for forecasting, model estimation, and inference; see Section \ref{sec: numerical} for details. Lastly, it is important to emphasize that $\alpha$ is not a hyperparameter; instead, it is a parameter that is determined immediately by the model choice, as detailed in the equation (\ref{params}). Also, the hyperparameter $\gamma$ only matters when we consider the crossing statistics-based model and in fact, does not play a role in predictability- and portmanteau criterion-based models.
}

{\color{black}Despite the extensive literature on solving close problems in portfolio optimization \cite{coelho2023performance,hooshmand2023model,jing2022bi, sehgal2021robust},} the only approach for handling (\ref{pr: p-original}) suggests to apply the following semidefinite (SDP) relaxation \cite{cuturi2016mean}:
\begin{equation}\label{pr: p-original-sdp-relaxation} \tag{\mbox{$SDP-P$}}
\begin{aligned}
\min_{X\in \mathbb{R}^{n\times n}} \quad & \ \alpha \tr(A_0X)+\gamma \sum_{i=2}^q \tr(A_iX)^2+\beta \|X\|_1 
\\
\textrm{subject to} \quad & 
 \tr(A_0X)\ge \phi, \quad
  \tr(X)=1,\quad \text{and} \quad    X\succeq 0,   
\end{aligned}
\end{equation}
with $\beta>0$ and $\|X\|_1:=\sum_{i,j}|X_{i,j}|$.
{\color{black}
However, this method exhibits several drawbacks: (i) optimal values of (\ref{pr: p-original}) and (\ref{pr: p-original-sdp-relaxation}) may differ,  (ii) a solution of (\ref{pr: p-original-sdp-relaxation}) is not necessarily rank-one, and (iii) even if it is rank-one, a common approach is to apply sparse principal component analysis to it but, the qualitative properties of a solution of the sparse component analysis is not well-understood with respect to the original problem. Therefore, the effectiveness of applying an SDP relaxation to (\ref{pr: p-original}) lacks rigorous theoretical support.}
{\color{black} The main reason is that the presence of $\|X\|_1$ in the objective function of (\ref{pr: p-original-sdp-relaxation}) breaks down the standard techniques that result in a relationship between a QCQP (without a cardinality constraint) with its original SDP relaxation. More specifically,  there are no known theoretical results that specify the relationship between the optimal values and points of a solution of a QCQP that has a cordiality constraint with its SDP relaxation. Further, there exists a body of literature addressing \textit{general} cardinality problems through other relaxation techniques like mixed-integer nonlinear programming \cite{bertsimas2009algorithm}, and continuous relaxations of mixed-integer programs \cite{burdakov2016mathematical} but such methods are less promising because they are not using the interesting structure of (P) in their design.}

Therefore, 
we propose an efficient tailored algorithm that
{\color{black} exploits the structure of (P) and}
utilizes a penalty decomposition method to directly solve (\ref{pr: p-original}) which obtains a KKT point of (\ref{pr: p-original}). 
{\color{black} To the best of our knowledge, the proposed algorithm is the first specialized approach specifically designed to directly tackle the problem (\ref{pr: p-original}) based on the portmanteau criterion and crossing statistics proxies.}
This unifying algorithm not only is applicable to all three proxies mentioned above but also works for nonconvex objective functions.
In the proposed methodology, we use the variable splitting technique to simplify our complex problem (\ref{pr: p-original}) that involves highly nonlinearly coupled constraints. Generally, the splitting approach involves introducing additional variables to the problem, which helps loosen the coupling between the original variables \cite{zeng2019global}. This transformation enables the problem to be tackled more efficiently. After splitting and penalizing over the related equality constraints, a sequence of penalized subproblems is in hand that could be approximately solved.  Each subproblem is efficiently tackled by a block coordinate algorithm that outputs a saddle point of the corresponding subproblem. Each step of the block coordinate algorithm either has a closed-form solution or is shown to be tractable; see problems (\ref{pr: px}), (\ref{pr: py}), and (\ref{pr: pz}), respectively, in 
Section \ref{sec: method}. 
In addition,  the sequence of objective values in the block coordinate algorithm is either strictly decreasing or reaches equality at a finite step, which also yields a saddle point. We establish that our tailored penalty decomposition algorithm finds a KKT  point of the original problem (\ref{pr: p-original}). Moreover, in the nonconvex case, the analysis of the algorithm can be simply extended to demonstrate that our algorithm is successful again, if the starting penalty parameter is larger than a finite positive bound and if the objective function has a bounded level set; see 
{\color{black} Theorem \ref{thm: nonconvex convergence} in Section \ref{sec: convergence} and its proof \ref{thm: proof-nonconvex} in the Appendix.}

{\color{black}
We mention that \cite{mousavi2022penalty} introduces a penalty decomposition method that exclusively addresses the sparse, volatile, and mean-reverting problem associated with the predictability proxy, specifically for when $\alpha=1$ and $\gamma=0$ in (\ref{pr: p-original}). In this work, we extend this methodology to efficiently tackle problems arising from other proxies as well.
}
We emphasize that, to the best of our knowledge, the portmanteau criterion- and crossing statistics-based problems with volatility and sparsity constraints have not been previously considered in the literature, even though they could produce more profitable portfolios as more information is leveraged by them. We numerically demonstrate that this is indeed the case for crossing statistics, which highlights the importance of designing an effective algorithm for this intractable problem; see Section \ref{sec: numerical}. Further, our empirical results reveal the quadratic term in the objective function is crucial in capturing mean reversion.

We demonstrate the effectiveness of the algorithm by conducting a numerical comparison with the SDP relaxation approach using the S\&P 500 dataset. Our method outperforms the SDP relaxation in terms of practical performance measures. Further, we provide a comparison of portfolio performance generated from predictability, portmanteau criterion, and crossing statistics proxies, when our algorithm is applied. The results show that after tuning the parameters $q$ and $\gamma$, the crossing statistics-based portfolio not only effectively captures the mean-reversion property, but also achieves superior performance in terms of Sharpe ratio and cumulative return and loss, compared to predictability and portmanteau criterion-based portfolios; see Section \ref{sec: numerical}. This is consistent with the fact that its corresponding problem exploits more information from the data in this case. However, the portmanteau criterion-based portfolio practically may not be successful in retaining mean reversion, underscoring the significance of the quadratic term in the objective function for capturing this desirable property.

\textbf{Organization.}
The remainder of the paper is organized as follows. Section \ref{sec: method} presents the proposed penalty decomposition algorithm for directly solving (\ref{pr: p-original}) together with all technical issues. Extensive numerical experiments are provided in \ref{sec: numerical}, while Section \ref{sec: conclusion} draws some concluding remarks. The proofs of the convergence analysis are presented in {\color{black} the Appendix}.

\textbf{Notation.}
The complement of a set $S$ is denoted as $S^c$ and its cardinality as $|S|$. For a natural number $n$, let   $[n]:=\{1,2,\dots,n\}$. 
For $S=\{i_1,i_2\dots, i_{|S|}\} \subseteq [n]$ and $x\in \mathbb R^n, x_S\in \mathbb R^{n}$ is the coordinate projection of $x$ with respect to indices in $S$, that is, $(x_S)_i=x_i$ for $i\in S$ and $(x_S)_i=0$ for  $i\in S^c$. 
By $\|.\|$, we mean the standard Euclidean norm. We show $A$ is positive semidefinite or definite by $A \succeq 0$ and $A\succ 0$, respectively.  
{\color{black} We denote the smallest and largest eigenvalues of $A$ with $\lambda_{\min}(A)$, and $\lambda_{\max}(A)$, respectively.}

\section{A Penalty Decomposition Method for Statistical Proxy-based Portfolios with Sparsity and Volatility Constraints} \label{sec: method}

 An algorithm for directly minimizing (\ref{pr: p-original}) that obtains a KKT point is presented next, with all proofs delegated to the Appendix. The PPC (standing for predictability, portmanteau criterion, and crossing statistics) based algorithm leverages a tailored penalty decomposition method in which a sequence of penalty subproblems is approximately solved. For each penalty subproblem, a block coordinate descent (BCD) algorithm is utilized that produces a saddle point of the penalty subproblem. The limit point of a suitable subsequence of such saddle points is demonstrated to be a KKT  point of (\ref{pr: p-original}).

\subsection{Details of PPC and BCD Algorithms}
By introducing two new variables $y$ and $z$, we can equivalently reformulate (\ref{pr: p-original}) as follows:
\begin{equation} \label{pr: equivalent_p} \tag{$P^\prime$}
\begin{aligned}
\min_{x,y,z\in \mathbb{R}^n} \quad & 
\alpha x^TA_1x+\gamma \sum_{i=2}^q (z^TA_iz)(x^TA_ix)
\\ \textrm{subject to} \quad & 
 x^TA_0x\ge \phi, \quad y^Ty=1, \quad  \|y\|_0\le k, \quad \\
 \quad & 
  x-y=0,\quad \text{and} \quad     x-z=0.       \\
\end{aligned}
\end{equation}
The reason for such splitting is to effectively handle the highly nonlinearly coupled constraints of the original problem. Specifically, we aim for the subproblems of Algorithm \ref{algo: BCD} to be as simple as possible.
Let
\begin{align} \label{def: q_rho}
q_{\rho}(x,y,z):=  \alpha x^TA_1x+\gamma \sum_{i=2}^q (z^TA_iz)  (x^TA_ix) + \rho(\|x-y\|_2^2+\|x-z\|_2^2),
\end{align}
and 
\begin{align*} 
\mathcal X:=\{x\in \mathbb{R}^n \, | \, x^TA_0x\ge \phi\}, \qquad
\mathcal Y:=\{y\in \mathbb{R}^n \, | \,y^Ty=1\, \text{ and } \|y\|_0\le k\}, \qquad \text{ and } \qquad
\mathcal Z :=\mathbb R^n.
\end{align*}
By penalizing the last two constraints in (\ref{pr: equivalent_p}),  we tackle this problem by a sequence of penalty subproblems as follows:
\begin{equation}\label{pr: pxyz} \tag{\mbox{$P_{x,y,z}$}}
\begin{aligned}
\min_{x,y,z}   \quad & \ 
\qquad q_{\rho}(x,y,z)
\\
\textrm{subject to}  &  \quad
  x\in \mathcal{X}, \quad y\in \mathcal{Y}, \quad \text{and} \quad z\in \mathcal{Z},   
\end{aligned}
\end{equation}
with $\rho$ going to infinity incrementally 
(such techniques have demonstrated their efficacy in the existing literature \cite{nazari2022penalty, mousavi2022penalty}).
The auxiliary variables $y$ and $z$ are introduced such that simpler subproblems can be obtained in Algorithm \ref{algo: BCD}. We establish in the sequel that this method  efficiently  finds a saddle point $(x_{*},y_{*},z_{*})$ of (\ref{pr: pxyz}), which means, $x_{*}\in \Argmin_{x\in \mathcal X} \,  q_{\rho}(x,y_*,z_*),$ 
 $ y_{*}\in\Argmin_{y\in \mathcal Y}  \, q_{\rho}(x_*,y,z_*),$ and 
 $z_{*}\in\Argmin_{z\in \mathcal Z} \, q_{\rho}(x_*,y_*,z).$

\begin{algorithm}
\caption{BCD Algorithm for Solving (\ref{pr: pxyz}) }
\begin{algorithmic}[1]
\label{algo: BCD} 
\STATE Input: Select arbitrary $y_{0}\in \mathcal Y$ and $z_{0}\in \mathcal Z$. 
\STATE Set $l=0$.
\STATE 
Solve $x_{l+1}\in\Argmin_{x\in \mathcal X} \ q_{\rho}(x,y_l,z_{l}).$
\STATE 
Solve $y_{l+1}= \Argmin_{y\in \mathcal Y} \ q_{\rho}(x_{l+1},y,z_{l}).$
\STATE 
Solve $z_{l+1}= \Argmin_{z\in\mathcal Z}\ q_{\rho}(x_{l+1},y_{l+1},z)$.
\STATE $l \leftarrow l+1$ and go to step (3).

\end{algorithmic} 
\end{algorithm}
Next, we discuss how to efficiently solve the restricted subproblems in Algorithm \ref{algo: BCD}.
\newline
$\bullet
\ \Argmin_{x\in\mathcal X}\ q_{\rho}(x,y,z)$: This subproblem becomes the following nonconvex quadratic optimization problem: 
\begin{equation} \label{pr: px} \tag{\mbox{$P_x$}}
\begin{aligned}
\min_{x\in \mathbb{R}^n}\quad & \ 
\alpha x^TA_1x+\gamma \sum_{i=2}^q \left(z^TA_iz\right)  \left(x^TA_ix\right)+\rho \left(\|x-y\|^2+\|x-z\|^2 \right)& 
\\ \textrm{subject to} & \quad \quad 
 x^TA_0x\ge \phi.
\end{aligned}
\end{equation}
Clearly, this problem has a solution $x_*$ and its KKT  conditions are as follows:
\begin{equation} \label{eqn: KKT_Px}
\alpha A_1x_* +\gamma \sum_{i=2}^q \left(z^TA_iz\right)\,  A_ix_*+\rho (2x_*-y-z)-\lambda A_0x_* =0, \quad
\mbox{ and } \quad 0\le \lambda \perp x_*^TA_0x_*-\phi \ge  0.
\end{equation}
In general, this nonconvex program does not have a closed-form solution, nevertheless, we will discuss how to find its \textit{global} minimizer by exploiting its SDP relaxation.
Recall that the  SDP  relaxation of a general quadratic program with exactly one general quadratic constraint obtains the same optimal value provided that it is strictly feasible \cite{boyd2004convex}. It is easy to see (\ref{pr: px}) is strictly feasible because  $A\succ 0$. Hence, we can  find the optimal value of this nonconvex problem exactly by solving its convex  SDP  relaxation given next:
\begin{equation} \label{pr: sdp_px}
\tag{\mbox{$SDP-P_x$}}
\begin{aligned}
\min_{X\in \mathbb R^{(n+1)\times (n+1)}}  \quad &
\tr\left(
H_1 X
\right)
\\
\textrm{subject to} \quad
 \tr&\left(H_2
 X\right)\ge 0, \quad X_{11}=1,
\quad \text{and} \quad
X\succeq 0,             \\
\end{aligned}
\end{equation}
where
$$
H_1=\begin{bmatrix} \rho\left( \|y\|^2+\|z\|^2\right) & -\rho (y+z)^T\\ -\rho (y+z) &  \alpha A_1+\gamma \sum_{i=2}^q(z^TA_iz)A_i+2\rho I \end{bmatrix}
\qquad \text{ and } \qquad 
H_2=\begin{bmatrix} -\phi & 0\\ 0& A_0 \end{bmatrix}.
$$
Its dual problem is given by
\begin{equation*}
\begin{aligned}
\max_{w_1 \in \mathbb R, w_2\in \mathbb R, Z\in \mathbb R^{(n+1)\times (n+1)}
} & w_2 
\\
\textrm{subject to} \qquad \quad &
w_1 H_2+w_2\begin{bmatrix} 1 & 0\\ 0& 0\end{bmatrix}+Z=H_1, \\
& w_1\ge 0, \quad \text{and} \quad
Z\succeq 0.           \\
\end{aligned}
\end{equation*}
We first show that both problems are strictly feasible. Clearly, since $A_0\succ 0$, we see $X=\begin{bmatrix} 1 & 0\\ 0& t I\end{bmatrix}$ with $t>\phi/\tr(A_0)>0$ is a strictly feasible point of the primal problem. To see that the dual problem is strictly feasible, it is enough to show that there exists a positive $w_1$ and a real $w_2$ such that $Z\succ 0.$
This block matrix is positive definite, if and only if (i) $\alpha A_1+\gamma \sum_{i=2}^q(z^TA_iz)A_i+2\rho I-w_1A_0\succ 0$ and (ii) $\rho (\|y\|^2+\|z\|^2)+w_1\phi-w_2 -\rho^2(y+z)^T( \alpha A_1+\gamma \sum_{i=2}^q(z^TA_iz)A_i+2\rho I-w_1A_0)^{-1}(y+z)>0$. To guarantee inequality (i),  it is enough to select $w_1=\epsilon>0$ small enough such that 
\begin{equation} \label{eqn: rho-px}
    \lambda_{\min}( \alpha A_1+\gamma \sum_{i=2}^q(z^TA_iz)A_i)+2\rho >\epsilon \lambda_{\max}(A_0),
\end{equation}
which is doable because $\lambda_{\min}(\alpha A_1+\gamma \sum_{i=2}^q(z^TA_iz)A_i))\ge 0$ in view of $A_i\succeq 0$ for $i\in [q], \gamma\ge0,$ and $\rho>0$. Inequality (ii) can be easily guaranteed by selecting a sufficiently large negative number for $w_2$.

Given that both the primal and the dual problems have strictly feasible solutions, their solutions have the same optimal value. Let $(w_1^*, w_2^*, Z^*)$ and $X^*$ be their optimal solutions. If $X^*$ of (\ref{pr: sdp_px}) has rank one, then the solution for (\ref{pr: px}) is trivially obtained. If not, using the procedure in \cite[Lemma 2.2]{ye2003new}, we can use the rank-one decomposition  $X^*=\sum_{i=1}^r u_iu_i^T$ with $r=\rank(X^*)$ and $ 0 \ne u_i \in \mathbb R^{n+1}$ for all $i\in [r]$ such that $ru_i^TH_2u_i=\tr\left(H_2X^*\right)\ge 0$ for all $i\in [r]$.
Since $X^*_{11}=1$, there exists $j\in [r]$ such that $u_j=[\alpha;u] \in \mathbb R^{n+1}$ with $\alpha\ne 0$. Further, the KKT  conditions imply that $0=\tr\left(X^*Z^*\right)=\tr\left(\sum_{i=1}^ru_iu_i^TZ^*\right)=\sum_{i=1}^r\tr\left(u^T_iZ^*u_i\right)=\sum_{i=1}^ru^T_iZ^*u_i$ and since $Z^*\succeq 0$, we have $u_j^TZ^*u_j=0$. Thus, $u_ju_j^T$ and $(w_1^*, w_2^*, Z^*)$  satisfy the KKT  conditions and consequently, $x_*:={u_j}/\alpha$ yields a \textit{global} solution to (\ref{pr: px}). Thus, $x_*$ indeed satisfies (\ref{eqn: KKT_Px}).
\newline
\newline
$\bullet
\ \Argmin_{y\in\mathcal Y}\ q_{\rho}(x,y,z)$: This subproblem is as follows:
\begin{equation} 
\label{pr: py} \tag{\mbox{$P_y$}}
\min_{y\in \mathbb{R}^n} \ \|x-y\|_2^2 \qquad \mbox{subject to}  \qquad  y^Ty=1, \quad \textrm{and} \quad \|y\|_0\le k.
\end{equation}
Thus, due to Lemma 3.2 in \cite{mousavi2022penalty},   the closed-form solution of this problem is given by
\begin{equation} \label{eqn: y_star_Py}
y_*=\mathcal T_k\left(x\right),
\end{equation}
where $\mathcal{T}_k$ is defined in (\ref{eqn: sparsifying}). 
\newline
\newline
$\bullet
\ \Argmin_{z\in\mathcal Z}\ q_{\rho}(x,y,z)$: This subproblem becomes the following:
\begin{equation} 
\label{pr: pz} \tag{\mbox{$P_z$}}
\min_{z\in \mathbb{R}^n} \ \gamma \sum_{i=2}^q \left(x^TA_ix\right)\left(z^TA_iz\right)+\rho\|x-z\|^2,
\end{equation}
which has the closed-form solution of 
\begin{equation} \label{eqn: z_star_BCD}
  z_*=\rho\left(\gamma
\sum_{i=2}^q \left(x^TA_ix\right)A_i
+\rho I\right)^{-1}x.  
\end{equation}

Next, we develop our proposed PPC Algorithm \ref{algo: PPC} that starts from an initial positive penalty parameter and smoothly enlarges it to numerically go to infinity. For each individual penalty parameter, Algorithm \ref{algo: BCD} is utilized to solve the corresponding subproblem. We assume that (\ref{pr: p-original}) is feasible and a feasible point $x^{\text{feas}}$ is available.  To find such a point, one can solve the following sparse 
{\color{black}
principal component analysis} problem:
\begin{equation*}
\max_{x\in \mathbb{R}^n} \ x^TA_0x \qquad \mbox{subject to}  \qquad  x^Tx=1 \quad \text{and}\quad \|x\|_0\le k.
\end{equation*}
A stationary point of this problem can be obtained by the power method $x_{l}:= \mathcal{T}_k(A_0x_{l-1})$ for $l\in \mathbb N$ (with $u_0\in \mathbb R^n$ arbitrary) suggested in \cite{ruszczynski2011nonlinear}. We can run this cheap method starting from different $x_0$'s until the achieved stationary point satisfies the volatility constraint. 
The sparsifying operator $\mathcal{T}_k$ is defined as  
\begin{equation}  \label{eqn: sparsifying}
\mathcal T_k(x):=\frac{x_\Lcal}{\left\|x_\Lcal\right\|_2},
\end{equation}
where $\Lcal$ is the index set containing the $k$  largest components of $x$ in absolute value. 
To provide this algorithm and its analysis later, we let:
\begin{align} \label{def: upsilon}
\Upsilon  \ge\max\{f(x^{\text{feas}}), \min_{x\in \mathcal{X}} q_{\rho^{(0)}}(x,y^{(0)}_0,z^{(0)}_0) \}> 0 
\qquad \text{ and } \qquad 
X_\Upsilon :=\{x \in \mathbb R^n\, | \, f(x)\le \Upsilon\},
\end{align}
{\color{black}{where $f(x)$ is defined in the objective function of (\ref{pr: p-original}).}}

\begin{algorithm}
\caption{PPC Penalty Decomposition Algorithm for Solving (\ref{pr: p-original}) }
\begin{algorithmic}[1]
\label{algo: PPC}
\STATE Inputs:  $r>1, \rho^{(0)}>0$, and  $y^{(0)}_0\in \mathcal Y$ and $z^{(0)}_0\in \mathcal Z$ with $\|z_0^{(0)}\|\le 1$. 
\STATE Set $j=0$.
\REPEAT
 \STATE Set $l=0$.
 \REPEAT
\STATE Solve $x^{(j)}_{l+1}\in \Argmin_{x\in \mathcal X} \ q_{\rho}(x,y^{(j)}_l,z^{(j)}_{l}).$
\STATE Solve
$y^{(j)}_{l+1}=\Argmin_{y\in \mathcal Y} \ q_{\rho}(x^{(j)}_{l+1},y,z^{(j)}_{l}).$
 \STATE Solve  $z^{(j)}_{l+1}=\Argmin_{z\in\mathcal Z}\ q_{\rho}(x^{(j)}_{l+1},y^{(j)}_{l+1},z)$.
 \STATE Set $l \leftarrow l+1$.
 \UNTIL{stopping criterion (\ref{BCD-practical-stopping-criterion}) is met.} \STATE Set
 $(x^{(j)},y^{(j)} ,z^{(j)}):= (x^{(j)}_{l}, y^{(j)}_{l},z^{(j)}_{l})$.

\STATE  Set $\rho^{(j+1)} = r\cdot \rho^{(j)}$.

\STATE 
If $\min_{x\in \mathcal{X}} q_{\rho^{(j+1)}}(x,y^{(j)},z^{(j)})> \Upsilon$, then $y_0^{(j+1)}=x^{\text{feas}}$ and $z_0^{(j+1)}=x^{\text{feas}}$. Otherwise,
$y^{(j+1)}_0 = y^{(j)}$ and $z^{(j+1)}_0 = z^{(j)}$.
\STATE Set $j \leftarrow j+1$.
\UNTIL{stopping criterion (\ref{outer loop stopping criteria}) is met}. 
\end{algorithmic}
\end{algorithm}
{\color{black}Note that Algorithm \ref{algo: PPC} consists of two nested loops,w}e stop the inner loop,
{\color{black} lines 5 to 10,}
if (after dropping $(j)$ for simplicity):
\begin{equation}\label{BCD-practical-stopping-criterion}
\max
\left\{
\frac{\|x_l-x_{l-1}\|_\infty}{\max \left(\|x_l\|_\infty,1 \right)},
\frac{\|y_l-y_{l-1}\|_\infty}{\max \left(\|y_l\|_\infty,1 \right)},
\frac{\|z_l-z_{l-1}\|_\infty}{\max \left(\|z_l\|_\infty,1 \right)}
\right\}
\le
\epsilon_I,
\end{equation}
and the outer loop, 
{\color{black}lines 3 to 15,}
is stopped  when the following convergence criterion is met:
\begin{equation} \label{outer loop stopping criteria}
\|x^{(j)}-y^{(j)}\|_{\infty}+\|x^{(j)}-z^{(j)}\|_{\infty}
\le
\epsilon_O.
\end{equation}

\subsection{Convergence Analysis}
We first investigate Algorithm \ref{algo: BCD} and then establish the convergence of Algorithm \ref{algo: PPC}.

\textbf{Analysis of Algorithm \ref{algo: BCD}.}
We analyze a sequence $\{(x_l,y_l,z_l)\}$ generated by Algorithm \ref{algo: BCD} and provide a tailored proof that any such sequence obtains a saddle point of (\ref{pr: pxyz}). This justifies the use of Algorithm \ref{algo: BCD} for this nonconvex problem.
\begin{lemma} \label{lem: BCD boundedness}
Let $\alpha\in \{0,1\},\gamma\ge 0, \rho>0, A_0\succ 0,$ and $A_i\succeq 0$ for $i\in [q]$.  Consider the iterates of Algorithm \ref{algo: BCD}, that is,   $x_{l}\in \Argmin_{x\in \mathcal X} \ q_{\rho}(x,y_{l-1},z_{l-1})$, $y_{l}\in \Argmin_{y\in \mathcal Y} \ q_{\rho}(x_{l},y,z_{l-1})$
and  further, $z_{l}\in \Argmin_{z\in\mathcal Z}\ q_{\rho}(x_{l},y_{l},z),
$ we get
\begin{equation} \label{eqn: xyz_ bounded above}
    \max\{\|x_{l}\|,\|y_l\|, \|z_{l}\|\}\le 
\max\{\sqrt{\phi/{\lambda_{\min}(A_0)}},
      1\}.
\end{equation}
\end{lemma}
This lemma establishes that the sequence created by Algorithm \ref{algo: BCD} is bounded. Therefore, every sequence generated by this algorithm possesses a point of accumulation. The subsequent theorem establishes that every accumulation point is unquestionably a saddle point of (\ref{pr: pxyz}). Additionally, the sequence $(q_\rho(x_l, y_l,z_l))$ is either strictly decreasing or two consecutive terms produce the same value, resulting in a saddle point. Essentially, Algorithm \ref{algo: BCD} produces a saddle point either in finite steps or in the limit.
\begin{theorem} \label{thm: BCD convergence}
Let $\{(x_l,y_l,z_l)\}$ be a sequence generated by Algorithm \ref{algo: BCD} for solving (\ref{pr: pxyz}). Suppose that $(x_*,y_*,z_*)$ is an accumulation point of this sequence, then $(x_*,y_*,z_*)$ is a saddle point of the nonconvex problem (\ref{pr: pxyz}). Moreover, $\left\{q_{\rho}(x_{l},y_{l},z_{l})\right\}$ is a non-increasing sequence. 
If $q_\rho(x_r, y_r,z_r) =  q_\rho(x_{r+1}, y_{r+1},z_{r+1})$ for some $r\in \mathbb N$, then $(x_r, y_r,z_r)$ is a saddle point of $(\ref{pr: pxyz})$.
\end{theorem}

\textbf{Analysis of Algorithm \ref{algo: PPC}.}
Suppose that $x^*$ is a local minimum of (\ref{pr: p-original}), then there exists an index set $\Lcal$ such that $|\Lcal|=k$ and $x^*_{\Lcal^c}=0$ such that $x^*$ is also a local minimizer of  the following problem:
\begin{equation*}
\begin{aligned}
\min_{x\in \mathbb{R}^n}\quad & \quad
 f(x)  
\\ \textrm{subject to} & \quad 
 x^TA_0x\ge \phi, \quad 
x^Tx=1,    \quad 
\text{and} \quad   
x_{\Lcal^c}=0,
\end{aligned}
\end{equation*}
{\color{black}{where $f(x)$ is defined in the objective function of (\ref{pr: p-original}).}}
Recall that Robinson's  constraint qualification conditions for a local minimizer $x^*$ for the above problem are as follows  \cite{mousavi2022penalty}:
\begin{flalign} \hfill
\label{eqn: robinson}
& (i)  \mbox{ If } (x^*)^T A_0 x^*> \phi, \text{ then }  \bigg\{  \begin{bmatrix} -2 d^T  A_0 x^*  - v \\ 2 d^T x^* \\ d_{\Lcal^c} \end{bmatrix} \  \Big | \   d \in \mathbb R^{n}, \  v \in \mathbb R \bigg\} = \mathbb R \times \mathbb R \times \mathbb R^{|\Lcal^c|}; \notag
\\
&(ii)  \mbox{ If }  (x^*)^T A_0 x^*= \phi, \text{ then } \bigg\{  \begin{bmatrix} -2 d^T  A_0 x^*- v \\ 2 d^T x^* \\ d_{\Lcal^c} \end{bmatrix} \  \Big | \   d \in \mathbb R^{n}, \  v \in \mathbb R_- \bigg\} = \mathbb R \times \mathbb R \times \mathbb R^{|\Lcal^c|}.
\end{flalign} 
It can be seen that $x^*_\Lcal \ne 0$, due to   $\|x^*\|_2=1$ and $x^*_{\Lcal^c}=0$, and thus Robinson's conditions for case (i) always hold. For case (ii), Robinson's conditions hold, if and only if $\{ (A_0 x^*)_\Lcal, x^*_\Lcal\}$ is linearly independent. It is easy to see this set is linearly independent almost always. This implies that except for a set of measure zero, Robinson's conditions are satisfied for (\ref{pr: p-original}), that is, such an assumption for $x^*$ is essentially the case in practice.

Under Robinson's conditions mentioned above, the KKT  conditions for a local minimizer $x^*$ of (\ref{pr: p-original}) are the existence of
{\color{black} Lagrangian multipliers}
$(\lambda,\mu,w)\in \mathbb R \times \mathbb R\times \mathbb R^n $ and $\mathcal{ L}\subseteq [n]$ with $|\mathcal{L}|=k$ such that the following holds:
\begin{equation}  \label{eqn: KKT _conditions for p}
\begin{aligned}
& \quad (\alpha A_1+2\gamma \sum_{i=2}^q(x^TA_ix)A_i)x-\lambda A_0x+\mu x+w=0,
\\
& x_{{\mathcal{L}}^c}=0, \, 0\le \lambda \perp x^TA_0x-\phi\ge 0, \, \|x\|_2=1, \text{ and } w_{\mathcal{L}}=0. 
\end{aligned}
\end{equation}
The next result states that an arbitrary sequence  $\{(x^{(j)}, y^{(j)},z^{(j)})\}$ generated by Algorithm \ref{algo: PPC} has a convergent subsequence, whose limit point is a KKT  point of (\ref{pr: p-original}).
\begin{theorem} \label{thm: PPC convergence}
Suppose that  $\alpha\in \{0,1\}, \gamma\ge 0,\phi >0, k\in [n], A_0\succ 0$  and  $A_i\succeq 0$ for $i\in [q]$.
Let  $\left\{\left(x^{(j)}, y^{(j)},z^{(j)}\right)\right\}$  be a sequence generated by Algorithm~\ref{algo: PPC} for solving (\ref{pr: p-original}). Then, the following hold:
 \begin{itemize}
  \item [(i)] $\left\{\left(x^{(j)}, y^{(j)},z^{(j)}\right)\right\}$ has a convergent subsequence whose accumulation point  $(x^*, y^*,z^*)$ satisfies $x^*=y^*=z^*$. Further, there exists an index subset $\Lcal\subseteq [n]$ with $|\Lcal|=k$ such that $x^*_{\Lcal^c}=0$.
  \item [(ii)] Suppose that Robinson's condition given in (\ref{eqn: robinson}) holds at $x^*$ with the index subset $\Lcal$ indicated above. Then, $x^*$ is a KKT  point satisfying (\ref{eqn: KKT _conditions for p}) for (\ref{pr: p-original}).
 \end{itemize}
\end{theorem}
\textbf{Extension to Nonconvex Objective Functions.} 
Note that the focus has been on constructing portfolios based on the predictability, portmanteau criterion, and crossing statistics proxies, for which $f(x)$ as defined in (\ref{def: upsilon}) is \textit{convex}. Nevertheless, we point out that from a technical viewpoint, one can employ Algorithm \ref{algo: PPC} to solve (\ref{pr: p-original}) even when the objective function is \textit{not} convex.  For this case, we require the boundedness of the level set $X_\Upsilon$ defined in (\ref{def: upsilon}) and $\rho^{(0)}$ to be large enough in this algorithm. 
\begin{theorem} \label{thm: nonconvex convergence} 
Suppose that  $\alpha\in \{0,1\}, \gamma\ge 0,\phi >0, A_0\succ 0, A_i$ be symmetric for $i\in [q],  k\in [n], \rho^{(0)}>  |\lambda_{\max}(A_0)-\alpha\lambda_{\min}(A_1)|$, and the level set $X_\Upsilon$ in (\ref{def: upsilon}) is bounded. All the steps in Algorithm \ref{algo: PPC} are well-defined and it successfully finds a KKT point of  (\ref{pr: p-original}).
\end{theorem}

\section{Numerical Experiments} \label{sec: numerical}
Here, we examine the performance of Algorithm \ref{algo: PPC} in finding sparse, volatile, and statistically proxy-based mean-reverting portfolios on real data coming from the US stock market Standard and Poor's (S$\&$P 500) Index. This data set is often used for showing the practical effectiveness of an algorithm or for comparing numerical performance between competing approaches. 

Recall that a statistical arbitrage strategy typically involves four main steps, namely constructing an appropriate asset pool, designing a mean-reverting portfolio, verifying the mean-reverting property through a unit-root test, and finally trading the favorable portfolio. The construction of an asset pool can be achieved using a method described in \cite{cuturi2016mean}, which utilizes the smallest eigenvalue of the covariance matrix. In this paper, we focus on the crucial second step of developing a mean-reverting portfolio by minimizing predictability, portmanteau criterion, or crossing statistics measures, while incorporating volatility and sparsity constraints to ensure the profitability of the portfolio in practice. The mean-reversion property can be verified using a unit-root test such as the Dickey-Fuller test \cite{dickey1979distribution}. A comprehensive discussion on the trading strategy of a mean-reverting portfolio in \cite{zhao2018optimal,zhao2018mean} is based on observed/normalized spread values. Normalizing the spread values makes the strategy less sensitive to absolute price levels. Table 1 in \cite{zhao2018mean} provides examples of how the strategy works for different trading positions and normalized spread values. The strategy involves closing long positions and taking short positions when the normalized spread value exceeds a threshold, and not taking any action when the normalized spread value is negative. It is a simple yet effective strategy, but careful testing and analysis are necessary before adopting it as a long-term trading strategy.

In our numerical experiments, various standard performance metrics are used
{\color{black} to measure the quality of portfolios, which in our case, are obtained from the limit point of Algorithm \ref{algo: PPC} and applying sparse PCA to the solution of (SDP-P). The first measure is} called cumulative profit and loss ($P\&L$) is used to measure the overall return of a mean-reverting portfolio within a single trading period, from $t_1$ to $t_2$, and is calculated as
$$
\text{Cumulative P\&L} (t_1,t_2)= \sum_{t=t_1}^{t_2}P\&L_t,
$$
where $P\&L_t= \color{black}{x}^Tr_t(t-t_o)-x^Tr_{t-1}(t-1-t_o)$ when a long position is opened, and $P\&L_t= \color{black}x^Tr_t(t-t_o)-x^Tr_{t-1}(t-1-t_o)$ when a short position is opened at time $t_o$. For a particular asset, the value of $r_t(\tau)=\frac{p_t-p_{t-\tau}}{p_{t-\tau}}\approx \ln (p_t)-\ln (p_{t-\tau})$, where $p_t$ represents the price of the asset at time $t$.  Table 1 in \cite{zhao2018mean} is used for this purpose with $d$ being the standard deviation of the portfolio.
The Return on Investment (ROI) is another metric used to evaluate the investment return of a mean-reverting portfolio. It is calculated as follows:
$\text{ROI}_t=\frac{P\&L_t}{\|x\|_1},$
where $P\&L_t$ is the profit and loss of the portfolio at time $t$, and $\|x\|_1$ is the L1-norm of the portfolio weights. The last metric used to evaluate the performance of a portfolio over a given period of time from $t_1$ to $t_2$ is the Sharpe ratio ($SR$), which is defined as 
$$
SR_{ROI}(t_1,t_2)=\mu_{ROI}/\sigma_{ROI},
$$
where $\mu_{ROI}=1/(t_2-t_1)\sum_{t=t_1}^{t_2} ROI_t$, and $\sigma^2_{ROI}=1/(t_2-t_1)\sum_{t=t_1}^{t_2} (ROI_t-\mu_{ROI})^2$. A portfolio with a higher $SR$ is considered more profitable.

To conduct the numerical experiments, we start by explaining how to select matrices $A_i$'s for $i\in\{0,1,\dots,q\}$. Given a multivariate stochastic process $\bf x=(\bf x_t)_{t\in \mathbb N}$ with values in $\mathbb R^n$. For a nonnegative integer $s$, the lag-$s$ \textit{empirical} autocovariance matrix of a sample path $x=(x_1,x_2,\dots, x_T)$ of $\bf x_t$ is defined as 
\begin{equation} \label{def: autocov}
\Gamma_s:= \frac{1}{T-s-1} \sum_{t=1}^{T-s}  \tilde{x_t}\tilde{x}_{T+s}^T \qquad \text{ with }  \qquad 
\tilde{x}_t:=x_t-\frac{1}{T} \sum_{t=1}^Tx_t.
\end{equation}
{\color{black}
In particular, we have 
\begin{equation}  \label{params}
 \left\{ 
\begin{array}{ll}
\alpha=1, \, \gamma=0, \ \text{and} \ A_0=\Gamma_0  & \quad \mbox{for the predictability-},
\\
\alpha=0,\, \gamma=1,\ \text{and} \ A_i=\Gamma_i; \ \forall i\in \{0,2,3,\dots,q\}
 & \quad \mbox{for the portmanteau criterion-, and}
\\
\alpha=1, \, \gamma>0,\ \text{and} \ A_i=\Gamma_i;\ \forall i\in \{0,1,2,3,\dots,q\}
 &  \quad \mbox{for the crossing statistics-based problems.}
\end{array} \right. \end{equation}}
We construct the asset pool by combining the pools suggested in optimization portfolio studies \cite{cuturi2016mean,zhao2018optimal,zhao2018mean,zhang2020sparse}, resulting in a pool of $n=30$ assets. The trading time period is from February 1st, 2012 to June 30th, 2014. 

We run Algorithm \ref{algo: PPC} with $r=\sqrt{10}, \epsilon_I= 10^{-3},$ defined in (\ref{BCD-practical-stopping-criterion}) and $\epsilon_O=10^{-3}$, defined in (\ref{outer loop stopping criteria}), and the volatility
threshold $\phi$ selected to be larger than thirty percent of the median variance of all assets in the pool (see \cite{cuturi2016mean}). In order to obtain reasonable results and compare the performance of the portmanteau criterion- and crossing statistics-based portfolios, we must decide upon a specific $q$. For this goal, we select the best $q$  among the candidate set  $\{2,3,4\}$ when we apply our algorithm for the case of the portmanteau criterion as follows. We measure the average Sharpe ratios for a small, medium, and large sparsity level $k=5,10,$ and $17$. Based on this procedure (results not reported due to space considerations), we get that $q=3$ produces better results than $q=2$, but the differences between $q=4$ and $q=3$ are rather marginal. Hence, in our experiments, we set $q=3$. 
{\color{black}
After fixing $q=3$, the next step involves selecting an appropriate $\gamma$. For this aim, we concentrate exclusively on portfolio optimization based on crossing statistics. We opt for a range of values within the interval $(0.0001, 1]$, incrementally adjusting $\gamma$ with a step size of 0.0009. While applying Algorithm 2 to the crossing statistic-based model across various sparsity levels, specifically $k=5, 10, and 17$, we evaluate the average Sharpe ratios. By this process, the best candidate is realized as $\gamma=0.001$. 
Furthermore, for (SDP-P), we need to determine $\beta$ as well.  To select this hyperparameter, we apply the methodology presented in \cite{cuturi2016mean} (that is, we solve this convex problem and apply sparse principal component analysis on its solution to introduce the corresponding portfolio), for all three cases reported in (\ref{params}) (with $q=3$, and $\gamma=0.001$) and values from 0.001 to 1.999, with 0.009 increments. This grid search revealed that the range between 0.8 and 1.2 showed better promise for $\beta$ values. We then narrowed down our selection within this range, using a 0.01 step. Throughout this iterative approach, we found that $\beta =1$ yielded the highest Sharpe ratio. As such, $\beta$ is selected as 1 in our comparisons.
}

\begin{figure}[H]
\centering  \includegraphics[width=0.49\textwidth,height=6.05cm]{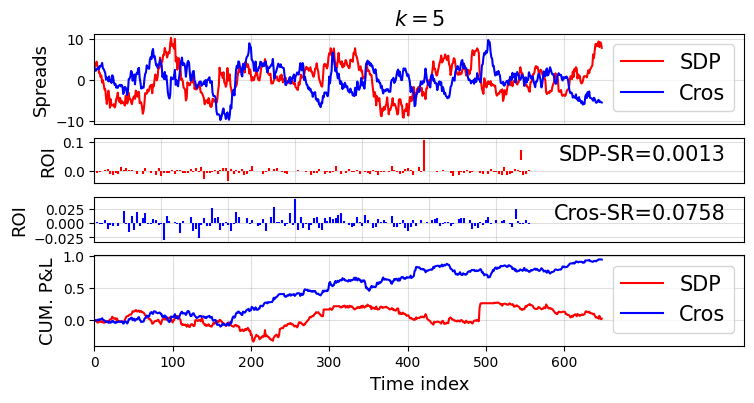}
\caption{Statistical proxy- and SDP-based portfolios for $k=5$.}
  \label{fig:k=5_time}
\end{figure}

\begin{figure}[H]
\centering
  \includegraphics[width=0.49\textwidth,height=6.05cm]{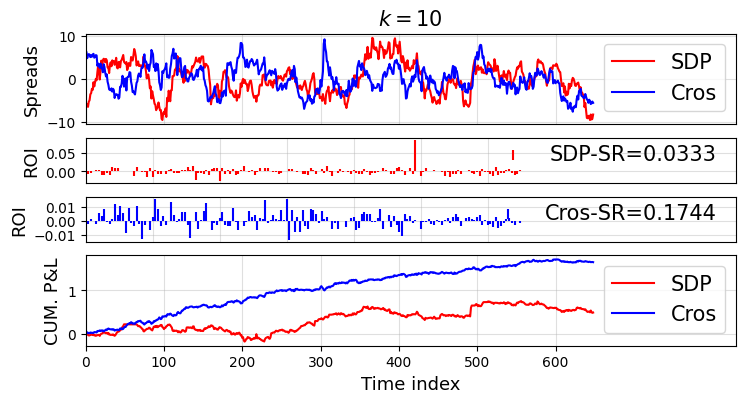}
\caption{Statistical proxy- and SDP-based portfolios for $k=10$.}
  \label{fig:k=10_time}
\end{figure}

\begin{figure}[H]
\centering
  \includegraphics[width=0.49\textwidth,height=6.05cm]{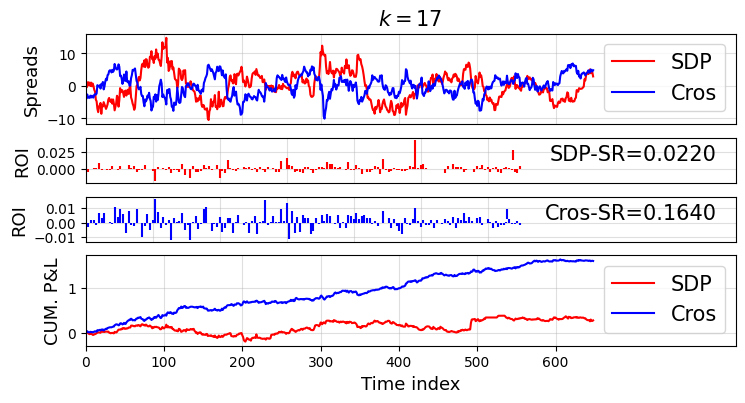}
\caption{Statistical proxy- and SDP-based portfolios for $k=17$.}
  \label{fig:k=17_time}
\end{figure}

{\color{black}
The spread pictures in Figures \ref{fig:k=5_time}, \ref{fig:k=10_time}, and \ref{fig:k=17_time} depict that the proposed method illustrates better mean reversion overall and the second and third pictures in these figures show that our method surpasses the SDP relaxation technique across both standard cumulative P\&L and Sharpe ratio assessments. This holds true for a range of sparsity levels, including small, medium, and large, as well as across all three portfolio models: predictability-based, portmanteau criterion-based, and crossing statistics-based portfolios.} Nevertheless, we should emphasize that the numerical performance of the SDP relaxation idea raises intriguing theoretical research questions that surely require thorough investigation. 
For example, under what conditions on the problem data in (\ref{pr: p-original}), its SDP relaxation (\ref{pr: p-original-sdp-relaxation}) has a rank-one solution. 
And, if available, how to use a rank-one solution of (\ref{pr: p-original-sdp-relaxation}) to construct a solution of (\ref{pr: p-original}).

\begin{figure}[H]
\centering
\includegraphics[width=0.49\textwidth,height=6.82cm]{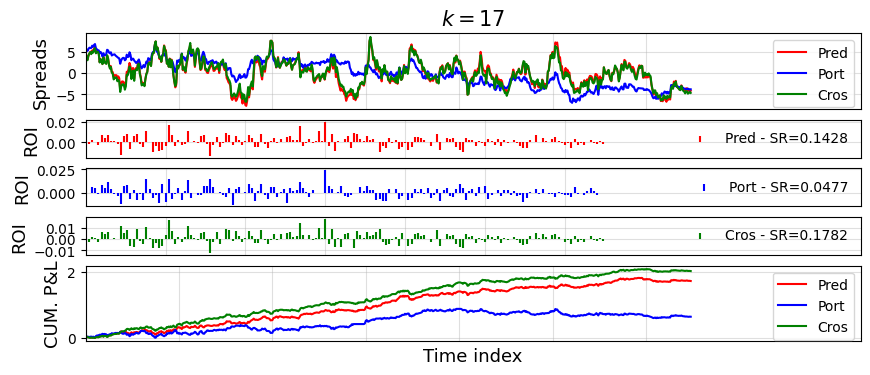}
\caption{Statistical proxy-based portfolios for $k=5,10,17$.}
\label{fig:proxy_portfolios}
\end{figure}

Next, we compare the performance of portfolios generated from the predictability, portmanteau criterion, and crossing statistics proxies, based on Algorithm \ref{algo: PPC}. Figure \ref{fig:proxy_portfolios} demonstrates that after careful tuning of the parameters $q$ and $\gamma$, the crossing statistics-based portfolio not only accurately captures the mean-reversion property, but also achieves superior performance in terms of Sharpe ratio and cumulative return and loss, surpassing the predictability and portmanteau criterion-based portfolios. This observation completely aligns with the idea that Algorithm \ref{algo: PPC} leverages more information from the data in this case. In addition, the portmanteau criterion proxy seems to fail in retaining mean reversion, underscoring the critical role of the quadratic term in the objective function for capturing this desirable property.

\begin{multicols}{2}
\begin{figure}[H]
\centering  \includegraphics[width=0.49\textwidth,height=6.05cm]{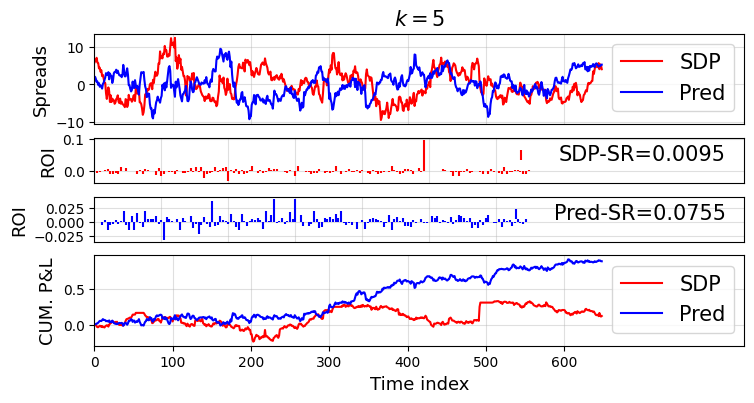}
\end{figure}
\begin{figure}[H]
\centering  \includegraphics[width=0.490\textwidth,height=6.05cm]{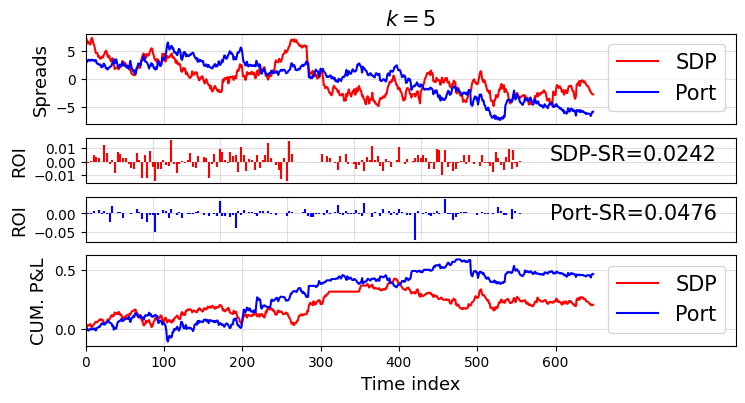}
\end{figure}

\begin{figure}[H]
\centering
  \includegraphics[width=0.49\textwidth,height=6.05cm]{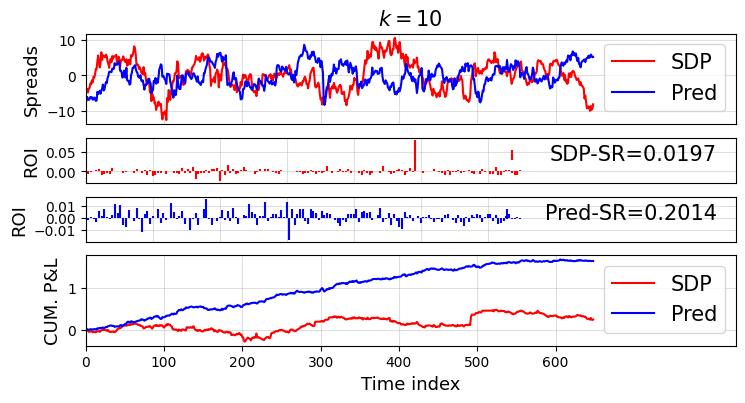}
\end{figure}
\begin{figure}[H]
\centering  \includegraphics[width=0.49\textwidth,height=6.05cm]{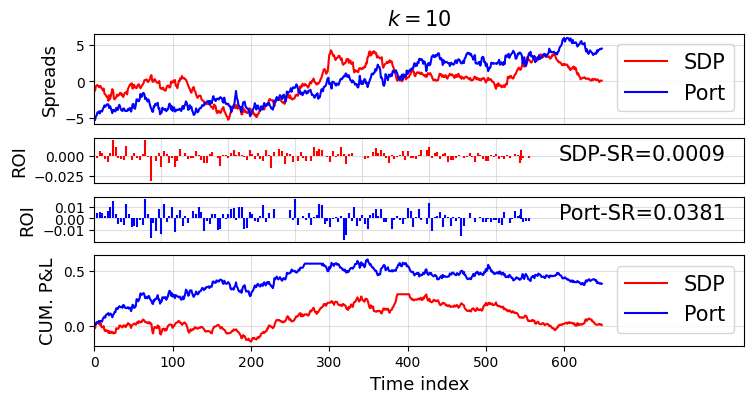}
\end{figure}

\begin{figure}[H]
\centering  \includegraphics[width=0.49\textwidth,height=6.05cm]{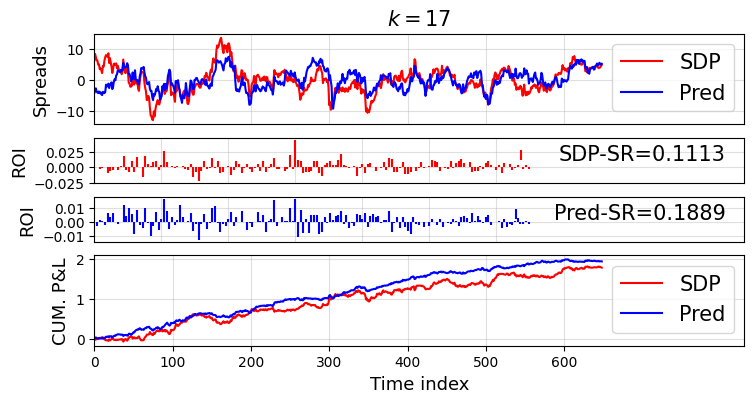}
\end{figure}
\begin{figure}[H]
\centering  \includegraphics[width=0.49\textwidth,height=6.05cm]{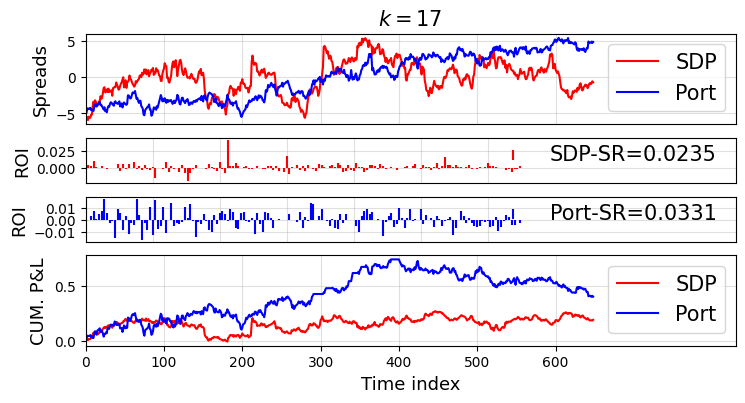}
\end{figure}

\begin{figure}[H]
\centering  \includegraphics[width=0.49\textwidth,height=6.82cm]{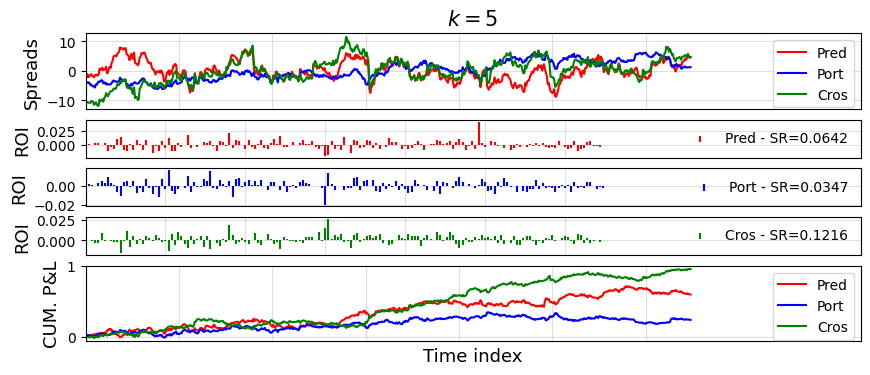}
\end{figure}
\begin{figure}[H]
\centering
  \includegraphics[width=0.49\textwidth,height=6.82cm]{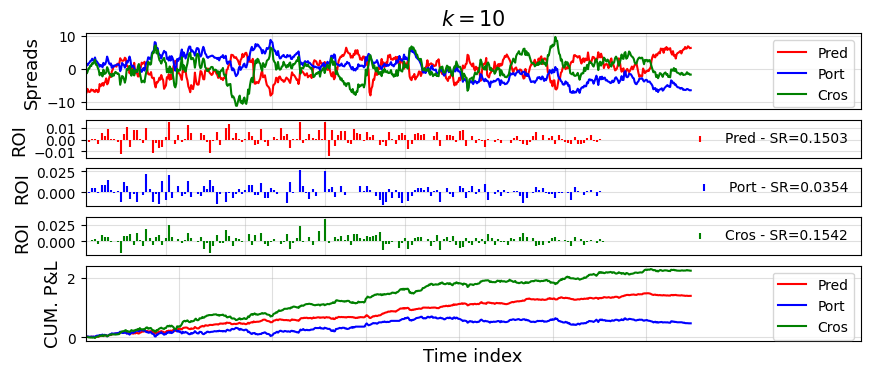}
\end{figure}

\end{multicols}

\section{Conclusion} \label{sec: conclusion}
The paper developed an algorithm that solves the underlying optimization problem of constructing portfolios of assets, having the mean-reverting property, and also being sparse and exhibiting reasonable volatility. Three main proxies for capturing mean-reversion are predictability, portmanteau criterion, and crossing statistics. A comprehensive optimization model was developed for this task.  We leverage the variable splitting technique to unfold highly coupled nonconvex constraints of this model and propose a tailored penalty decomposition method, that solves a sequence of penalty subproblems via an efficient block coordinate method. We establish that not only the proposed algorithm finds a KKT point when the objective function is convex,  but also its analysis can be extended to a nonconvex objective function (with partial modifications).  
We numerically examine the performance of the algorithm and show that it outperforms the SDP relaxation approach in terms of standard measures on an S\&P 500 data set. Further, we compare these statistical proxy-based portfolios and demonstrate that, after tuning hyperparameters,  the crossing statistics surrogate captures mean-reversion well and achieves better performance compared to predictability and portmanteau criterion-based portfolios. Our method has significant potential applications in finance, economics, and other related fields.

\bibliographystyle{plain}
 \bibliography{references}
\newpage
\section{Appendix} \label{sec: convergence}
All technical proofs are given next.
\subsection{Proof of Lemma \ref{lem: BCD boundedness}.}
\begin{proof} 
We drop the subscripts for simplicity. If $x^TA_0 x>\phi$, the Lagrangian multiplier in (\ref{pr: px}) must be zero such that we have 
$\left(\alpha A_1+\gamma\sum_{i=2}^q (z^TA_iz)A_i  + 2\rho I\right)  x ={\rho}(y+z
),$ 
and consequently, 
$x = \rho (\alpha A_1+\gamma\sum_{i=2}^q (z^TA_iz)A_i  + 2\rho I)^{-1}(y+z)$. 
This leads to
\begin{equation} \label{ineq: x_smaller_z}
\begin{aligned}
\|x\| 
&\le
\rho \| (\alpha A_1+\gamma\sum_{i=2}^q (z^TA_iz)A_i  + 2\rho I)^{-1}\|_2 \|y+z \|_2 
\\
&\le
\frac{\rho}{ \lambda_{\min}(\alpha A_1+\gamma\sum_{i=2}^q (z^TA_iz)A_i) + 2\rho} \|y+z \|_2 
\\& \le
\frac{1}{2} \|y+z \|_2  
\\
&\le
\frac{1}{2} (\|y\|+\|z\|)
\\
&=
\frac{1}{2}(1+\|z\|)  \qquad \text{as} \ \ y \ \ \text{solves} \ \  \ref{pr: py},
\end{aligned}
\end{equation} 
where we used that $ \lambda_{\min}(\alpha A_1+\gamma\sum_{i=2}^q (z^TA_iz)A_i)\ge 0$ in the third inequality.
Moreover, since $z$ solves (\ref{pr: pz}) and $\lambda_{\min}\left(\gamma \sum_{i=2}^q(x^TA_ix)A_i\right)\ge 0$,  the equation (\ref{eqn: z_star_BCD}) implies
\begin{equation} \label{ineq: z_smaller_x}
\begin{aligned}
\|z\| 
&\le
\rho \| (\gamma \sum_{i=2}^q(x^TA_ix)A_i+\rho I)^{-1}\| \|x\|
\\
&\le
 \frac{\rho}{\gamma \lambda_{\min}(\sum_{i=2}^q(x^TA_ix)A_i)+\rho}\|x\|
\\& \le
\|x\|.
\end{aligned}
\end{equation}
Through (\ref{ineq: x_smaller_z}), we see $\|x\|\le\frac{1}{2}(1+\|z\|)\le\frac{1}{2}(1+\|x\|)$; leading to $\|x\|\le 1$ in this case.  If $x^TA_0x=\phi$, we have $\|x\|\le \sqrt{ \phi/{\lambda_{\min}(A_0)}}$.  By (\ref{ineq: z_smaller_x}), we have $\|z\|\le \|x\|\le \max\{\sqrt{\phi/{\lambda_{\min}(A_0)}},1\}.$
\end{proof}

\subsection{Proof of Theorem \ref{thm: BCD convergence}.}
\begin{proof}
By observing definitions of $x_{l+1}, y_{l+1}$, and $z_{l+1}$ in steps 3-5 of Algorithm \ref{algo: BCD}, we get
\begin{eqnarray} \label{eqn: q_xyz-inequality}
  q_{\rho}(x_{l+1},y_{l+1},z_{l+1}) & \le & q_{\rho}(x_{l+1},y_{l+1},z); \qquad \forall z\in \mathcal Z, \notag
  \\ 
  q_{\rho}(x_{l+1},y_{l+1},z_l) &\le &  q_{\rho}(x_{l+1},y,z_{l}); \qquad \forall y\in \mathcal Y,\notag
  \\
  q_{\rho}(x_{l+1},y_{l},z_{l})  &\le &  q_{\rho}(x,y_{l},z_l); \qquad \forall x\in \mathcal X.
\end{eqnarray}
This simply leads to the following:
\begin{equation} \label{eqn: non-increasing_q}
\begin{aligned}
0  &\le q_{\rho}(x_{l+1},y_{l+1},z_{l+1}) \\
   & \le  q_{\rho}(x_{l+1},y_{l+1},z_{l})\\
   & \le q_{\rho}(x_{l+1},y_{l},z_{l})\\
   & \le q_{\rho}(x_{l},y_{l},z_{l}); \quad \forall l \in \mathbb N.
\end{aligned}
\end{equation}
Thus, $q_{\rho}(x_{l},y_{l},z_{l})$ is a  bounded  below and non-increasing sequence; implying that  $q_{\rho}(x_{l},y_{l},z_{l})$ is  convergent. 
From the other side, since $(x_*,y_*,z_*)$ is an accumulation point of 
 $\{(x_l,y_l,z_l)\}$, there exist a subsequence $L$ such that $\lim_{l\in L\to \infty} (x_l,y_l,z_l) = (x_*,y_*,z_*)$. The continuity of $q_{\rho}(x_{l},y_{l},z_{l})$ yields
\begin{equation*} 
\begin{aligned}
\lim_{l\to \infty}q_{\rho}(x_{l+1},y_{l+1},z_{l}) & = \lim_{l\to \infty} q_{\rho}(x_{l+1},y_{l},z_{l}) \\
   & =  \lim_{l\to \infty} q_\rho(x_l,y_l,z_l)\\
   & = \lim_{l\in L\to \infty} q_\rho(x_l,y_l,z_l) \\
   & = q_\rho(x_*,y_*,z_*).
\end{aligned}
\end{equation*}
By the continuity of  $q_{\rho}(x_{l},y_{l},z_{l})$ and taking the limit of both sides of (\ref{eqn: q_xyz-inequality}) as $l\in L\to \infty$ , we have 
\begin{eqnarray*} 
  q_{\rho}(x_{*},y_{*},z_{*}) & \le & q_{\rho}(x,y_{*},z_*); \qquad \forall x\in \mathcal X,\notag
  \\ 
  q_{\rho}(x_{*},y_{*},z_{*}) &\le &  q_{\rho}(x_{*},y,z_{*}); \qquad \forall y\in \mathcal Y,\notag
  \\
  q_{\rho}(x_{*},y_{*},z_{*})  &\le &  q_{\rho}(x_{*},y_{*},z); \qquad \forall z\in \mathcal Z.
\end{eqnarray*}
Further, it is clear from  (\ref{eqn: non-increasing_q}) that $\left\{q_{\rho}(x_{l},y_{l},z_{l})\right\}$ is non-increasing. 
\newline
Now suppose $q_\rho(x_l, y_l,z_l) =  q_\rho(x_{l+1}, y_{l+1},z_{l+1})$ for some $l\in \mathbb N$. Then because of the last inequality in (\ref{eqn: non-increasing_q}), we have $q_{\rho}(x_{l+1},y_{l},z_{l})=
q_{\rho}(x_{l},y_{l},z_{l}).$ 
Since $q_{\rho}(x_{l+1},z_l,y_l)=\min_{x\in \mathcal X} q_{\rho}(x,y_l,z_l)$ and $x_l\in \mathcal{X}$, we see $x_l\in \Argmin_{x\in \mathcal X}  q_{\rho}(x,y_l,z_l)$. 
Further, if  $q_\rho(x_l, y_l,z_l) =  q_\rho(x_{l+1}, y_{l+1},z_{l+1})$ for some $l\in \mathbb N$,  using the third inequality in (\ref{eqn: non-increasing_q}) we get
$q_{\rho}(x_{l+1},y_{l+1},z_{l}) = q_{\rho}(x_{l+1},y_{l},z_{l})$. Since $q_{\rho}(x_{l+1},y_{l+1},z_l)=\min_{y\in \mathcal Y} q_{\rho}(x_{l+1},y,z_l)$ and $y_l\in \mathcal{Y}$, we see $y_l\in \Argmin_{y\in \mathcal Y}  q_{\rho}(x_{l+1},y,z_l)$. By the last inequality in (\ref{eqn: non-increasing_q}), we have $  q_{\rho}(x_{l+1},y_l,z_l)=  q_{\rho}(x_{l},y_l,z_l)$, so $y_l\in  \Argmin_{y\in \mathcal Y}  q_{\rho}(x_{l},y,z_l)$.
Recall that $z_l \in \mathcal Z$ satisfies $z_l \in \mbox{Argmin}_{z \in \mathcal Z} \ q_\rho(x_l, y_l,z)$ by definition. Hence, $(x_l, y_l,z_l)$ is a saddle point of (\ref{pr: pxyz}). 
\end{proof}

\subsection{Proof of Theorem \ref{thm: PPC convergence}.}
\begin{proof}
(i)
From (\ref{eqn: xyz_ bounded above}), we know that
$    \max\left\{\|x^{(j)}\|,\|y^{(j)}\|, \|z^{(j)}\|\right\}\le 
      \max\{\sqrt{\phi/{\lambda_{\min}(A_0)}},1\},$
for every $j\in \mathbb N$. Thus,   $\{(x^{(j)}, y^{(j)}, z^{(j)})\}$ is bounded and therefore, has a convergent subsequence. For our purposes, without loss of generality, we suppose that the sequence itself is convergent. Let  $\left(x^*, y^*,z^*\right)$ be its accumulation point. 
First, let us show that $x^*=y^*=z^*$. Since  $\alpha\in \{0,1\}, A_i\succeq 0$ for $i\in [q],$ and $ \gamma\ge 0$, we see 
$ \alpha {(x^{(j)})}^TA_1x^{(j)}+\gamma \sum_{i=2}^q ({z^{(j)}})^TA_iz^{(j)} ({x^{(j)}})^TA_ix^{(j)}  \ge 0$, and thus in view of (\ref{def: q_rho}), (\ref{eqn: non-increasing_q}), and step 13 of Algorithm \ref{algo: PPC}, we have
$\rho^{(j)}\left(\|x^{(j)}-y^{(j)}\|^2+\|x^{(j)}-z^{(j)}\|^2 \right) \le \Upsilon$.
The latter leads to
$$
\max\{\|x^{(j)}-y^{(j)}\|,\|x^{(j)}-z^{(j)}\|\}\le \sqrt{{\Upsilon}/{\rho^{(j)}}},
$$
and thus $\max\{\|x^{(j)}-y^{(j)}\|,\|x^{(j)}-z^{(j)}\|\}\to 0$ when $j\to \infty$ as $\rho^{(j)}\to \infty$; implying that $x^*=y^*=z^*$.

Next, let $\Ical^{(j)}\subseteq [n]$ be such that $|\Ical^{(j)}|=k$ and $(y_{(\Ical^j)^c})_i=0$ for every $j\in \mathbb N$ and $i\in (\Ical^j)^c$. Then, since $\{\Ical^{(j)}\}$ is a bounded sequence of indices, it has a convergent subsequence, which means that there exists an index subset $\Lcal \subseteq [n]$ with $|\Lcal|=k$ and a subsequence $\{(x^{(j_\ell)}, y^{(j_\ell)}, z^{(j_\ell)})\}, $ of the above convergent subsequence such that $\Ical^{(j_\ell)} = \Lcal$ for all large $j_\ell$'s. Therefore, since $x^*=y^*$ and $y^*_{\Lcal^c}=0$, we see $x^*_{\Lcal^c}=0$, which further implies $ \|x^*_\Lcal\|=1$.

(ii)  For each $j$, $(x^{(j)}, y^{(j)}, z^{(j)})$ is a saddle point of  (\ref{pr: pxyz}) with $\rho=\rho^{(j)}>0$ such that (\ref{eqn: KKT_Px}),(\ref{eqn: y_star_Py}) and (\ref{eqn: z_star_BCD}) give
\begin{equation} \label{eqn: KKT -pxyz-convergence}
\left\{ \begin{aligned}
& \alpha A_1x^{(j_{\ell})}+\gamma \sum_{i=2}^q (z^{(j_{\ell})})^TA_iz^{(j_{\ell})} A_ix^{(j_{\ell})}+\dots \\
& \qquad \dots + \rho (x^{(j_{\ell})}-y^{(j_{\ell})}+x^{(j_{\ell})}-z^{(j_{\ell})})-\lambda A_0x^{(j_{\ell})} =0,
\\ 
&
y^{(j_{\ell})}={(x^{(j_{\ell})})_{\Lcal}}/{\|{(x^{(j_{\ell})})_{\Lcal}}\|},
\\
 &
 \gamma \sum_{i=2}^q (x^{(j_{\ell})})^TA_ix^{(j_{\ell})}A_i
z^{(j_{\ell})}=\rho (x^{(j_{\ell})}-z^{(j_{\ell})}),
\\  
&
0\le \lambda \perp (x^{(j_{\ell})})^TA_0x^{(j_{\ell})}- \phi \ge 0, \\
&\quad  \|y^{(j_{\ell})}\|=1, \quad \mbox{and} \quad (y^{(j_{\ell})})_{\Lcal^c}=0.
\end{aligned} \right.
\end{equation}

By injecting the first two lines of (\ref{eqn: KKT -pxyz-convergence}) into the third, we get 
\begin{equation*}
    \begin{aligned}
& \alpha A_1x^{(j_{\ell})}+\gamma \sum_{i=2}^q [(z^{(j_{\ell})})^TA_iz^{(j_{\ell})}  A_ix^{(j_{\ell})}
+ (x^{(j_{\ell})})^TA_ix^{(j_{\ell})}A_iz^{(j_{\ell})}]+
\rho^{(j_{\ell})} (x^{(j_{\ell})}-y^{(j_{\ell})})
-\lambda^{(j_{\ell})} A_0x^{(j_{\ell})}=0.
\end{aligned}
\end{equation*}
By observing
$
{\|(x^{(j_{\ell})})_\Lcal\|}\left(y^{(j_{\ell})}\right)_\Lcal=\left(x^{(j_{\ell})}\right)_{\Lcal}$ 
that implies 
$\left(x^{(j_{\ell})}-y^{(j_{\ell})}\right)_\Lcal=\underbrace{( \|x^{(j_{\ell})}\|-1 )
}_{:=\mu^{(j_\ell)}}(y^{(j_{\ell})})_\Lcal, 
$
and letting 
\begin{equation*}
    \begin{aligned}
M^{(j_\ell)}:= \alpha A_1x^{(j_{\ell})}+\gamma \sum_{i=2}^q [(z^{(j_{\ell})})^TA_iz^{(j_{\ell})}  A_ix^{(j_{\ell})}+(x^{(j_{\ell})})^TA_ix^{(j_{\ell})}A_iz^{(j_{\ell})}],
\end{aligned}
\end{equation*} 
we get
\begin{equation*}
 \label{eqn: KKT -subsequence-equation-passing-limit}
M^{(j_\ell)} + \mu^{(j_\ell)}  \begin{bmatrix}
 (y^{(j_\ell)})_{\Lcal} \\ 0 \end{bmatrix} +\underbrace{\begin{bmatrix}  0 \\ \rho^{(j_\ell)}(  x^{(j_\ell)})_{\Lcal^c} \end{bmatrix}  }_{:=w^{(j_\ell)}}-\lambda A_0x^{(j_{\ell})}=0,
\end{equation*} 
where $(w^{(j_\ell)})_\Lcal=0$ for each $j_\ell$ and also, $0\le \lambda^{(j_{\ell})}\perp(x^{(j_{\ell})})^TA_0x^{(j_{\ell})}-\phi\ge 0,\ \|y^{(j_{\ell})}\|_2=1,$ and $(y^{(j_{\ell})})_{\Lcal^c}=0.$
We next prove that $\{(\lambda^{(j_\ell)},\mu^{(j_\ell)},w^{(j_\ell)})\}$ is bounded  under Robinson's condition on $x^*$. 
 Suppose not, consider the normalized sequence
\[
( \wt \lambda^{(j_\ell)}, \wt \mu^{(j_\ell)}, \wt w^{(j_\ell)}):= \frac{ (\lambda^{(j_\ell)}, \mu^{(j_\ell)}, w^{(j_\ell)})} { \| (\lambda^{(j_\ell)}, \mu^{(j_\ell)}, w^{(j_\ell)}) \|_2}, \qquad \forall \ {j_\ell}.
\]
Through boundedness, there exists a convergent subsequence of $( \wt \lambda^{(j_\ell)}, \wt \mu^{(j_\ell)}, \wt w^{(j_\ell)})$ whose limit is given by $(\wt \lambda_*, \wt \mu_*, \wt w^*)$ such that $\| (\wt \lambda_*, \wt \mu_*, \wt w^*) \|_2=1$, we obtain, in view of $x^*=y^*$, $y^*_{\Lcal^c}=0$ and the boundedness of $(M ^{(j_\ell)})$, and  passing the limits, that
\begin{equation} \label{eqn:limit_condition}
   -\wt\lambda_*  A_0 x^* + \wt \mu_* x^* + \wt w^* =0,
\end{equation}
where $\wt \lambda_* \ge 0$, $x^*_{\Lcal^c} =0$, and $\wt w^*_\Lcal=0$. Let us consider two cases: $(x^*)^T A_0 x^* = \phi$, and $(x^*)^T A_0 x^* > \phi$ as follows.
When  $(x^*)^T A_0 x^* = \phi$, by the Robinson's conditions at $x^*$, there exist a vector $d \in \mathbb R^n$ and a constant $v\in \mathbb R_-$ such that $-2 d^T A_0 x^* - v = -2\wt \lambda_*$, $2 d^T x^* = - 2\wt \mu_*$, and $d_{\Lcal^c} = - \wt w^*_{\Lcal^c}$. Since $d_{\Lcal^c} = - \wt w^*_{\Lcal^c}$ and $\wt w^*_\Lcal=0$, we see that $d^T \wt w^* = -\| \wt w^* \|^2_2$. Therefore,
$0 = -\wt \lambda_* d^T A_0 x^* + \wt \mu d^T x^* + d^T \wt w^* =  - (\wt \lambda_*)^2 + \frac{\wt \lambda_* v }{2} - (\wt \mu_*)^2 - \|\wt w^* \|^2_2 = - \|  (\wt \lambda_*, \wt \mu_*, \wt w^*) \|^2_2 +  \frac{\wt \lambda_* v }{2},$
which implies that $\|  (\wt \lambda_*, \wt \mu_*, \wt w^*) \|^2_2 =  \frac{\wt \lambda_* v }{2}$.
Since $\wt \lambda_* \ge 0$ and $v \le 0$, we have $\|  (\wt \lambda_*, \wt \mu_*, \wt w^*) \|^2_2=0$, which is a  contradiction. Therefore, the sequence $\big( (\lambda^{(j_\ell)}, \mu^{(j_\ell)}, w^{(j_\ell)}) \big)$ is bounded.
Now, suppose $(x^*)^T A_0 x^* > \phi$. In this case,  since $(x^{j_\ell})$ converges to $x^*$, we have  $(x^{j_\ell})^T A_0 x^{j_\ell} > \phi$ for all $j_\ell$ sufficiently large. Hence, $\lambda^{(j_\ell)} =0$ for all large $j_\ell$. This shows that $\wt \lambda_*=0$.
In view of the Robinson's condition at $x^*$ given in (\ref{eqn: robinson}), there exist a vector $d \in \mathbb R^n$ and a constant $v\in \mathbb R$ such that $-2 d^T A_0 x^* - v = -2\wt \lambda_*$, $2 d^T x^* = - 2\wt \mu_*$, and $d_{\Lcal^c} = - \wt w^*_{\Lcal^c}$.  Using these equations (\ref{eqn:limit_condition}) together with $\wt \lambda_*=0$, we can see that
$0 = -\wt \lambda_* d^T A_0 x^* + \wt \mu d^T x^* + d^T \wt w^* =- \|  (\wt \lambda_*, \wt \mu_*, \wt w^*) \|^2_2 +  \frac{\wt \lambda_* v }{2};$ implying that $\|  (\wt \lambda_*, \wt \mu_*, \wt w^*) \|^2_2$, which is again a  contradiction.

Hence, $\{(\lambda^{(j_\ell)},\mu^{(j_\ell)},w^{(j_\ell)})\}$ is bounded and has a convergent subsequence with the limit $(\lambda,\mu, w)$. thus, through passing limit in (\ref{eqn: KKT -subsequence-equation-passing-limit}) and applying the results of part (i),  we have:\\
$\alpha A_1x^*+2\gamma \sum_{i=2}^q(x^*)^TA_ix^*A_ix^*-\lambda A_0x^*+\mu x^*+w=0, \|x^*\|=1, \quad (x^*)_{\Lcal^c}=0, \quad \mbox{and} \quad w^*_\Lcal =0,
$
where $ 0\le \lambda \perp (x^*)^TA_0x^*-\phi\ge 0$. This means that $x^*$ satisfies the first-order optimality conditions of (\ref{pr: p-original}) given in (\ref{eqn: KKT _conditions for p}).
\end{proof}

\subsection{Proof of Theorem \ref{thm: nonconvex convergence}.}
\begin{proof} \label{thm: proof-nonconvex}
Note that only a few arguments given in the above proofs must be partially modified in order to make sure that all the subproblems are well-defined and that the convergence analysis similarly holds in the nonconvexity case. Therefore, we avoid repetition and instead, shortly discuss all the arguments that must be revisited.
Specifically, we must put conditions on $\rho$ as follows.
\newline
(1) $ 2\rho > \lambda_{\max}(A_0)-\lambda_{\min}( \alpha A_1+\gamma \sum_{i=2}^q(z^TA_iz)A_i)$, which is to guarantee the existence of a  strictly feasible point in the dual of the SDP relaxation of (\ref{pr: px}); simply inequality (\ref{eqn: rho-px}). (2) $\rho>-\gamma \lambda_{\min}( \sum_{i=2}^q (x^TA_ix)A_i)$, which is to have well-defined $z_*$  in (\ref{eqn: z_star_BCD}). And, 
(3) $ 2\rho> -\lambda_{\min}(\alpha A_1+\gamma\sum_{i=2}^q (z^TA_iz)A_i)$, which is for the inequality (\ref{ineq: x_smaller_z}). 
Since $A_0\succ 0$, (3) is satisfied if (1) holds. It is easy to show that conditions (1) and (2) are satisfied if $\rho >  \lambda_{\max}(A_0)-\alpha \lambda_{\min}(A_1).$ 
Recall that for two symmetric matrices $Q_1$ and $Q_2$, we   $x^TQ_1x\ge \lambda_{\min}(Q_1)\|x\|^2$ and for two symmetric matrices $\lambda_{\min}(Q_1+Q_2)\ge \lambda_{\min}(Q_1)+\lambda_{\min}(Q_2)$. Thus, we can show 
$\lambda_{\min}(\alpha A_1+\gamma\sum_{i=2}^q (z^TA_iz)A_i)\ge \alpha \lambda_{\min}(A_1) +\gamma \sum_{i=2}^q (z^TA_iz)\lambda_{\min}(A_i)
$
 $ \ge \alpha \lambda_{\min}(A_1) +\gamma \sum_{i=2}^q \|z\|^2\lambda^2_{\min}(A_i)
\ge \alpha \lambda_{\min}(A_1);
$ implying that $\lambda_{\max}(A_0)-\lambda_{\min}( \alpha A_1+\gamma \sum_{i=2}^q(z^TA_iz)A_i)\le \lambda_{\max}(A_0)-\alpha\lambda_{\min}(A_1)$. 
In the same manner 
$\lambda_{\min}( \sum_{i=2}^q (x^TA_ix)A_i)\ge \sum_{i=2}^q (x^TA_ix)\lambda_{\min}(A_i)\ge \sum_{i=2}^q \|x\|^2\lambda^2_{\min}(A_i)\ge 0$ such that we have
$-\gamma \lambda_{\min}( \sum_{i=2}^q (x^TA_ix)A_i)\le 0$. Therefore, conditions (1) and (2) are satisfied if $\rho> |\lambda_{\max}(A_0)-\alpha\lambda_{\min}(A_1)|.$
\newline
Moreover, for part (i) in Theorem \ref{thm: PPC convergence}, we can have 
$\alpha {(x^{(j)})}^TA_1x^{(j)}+\gamma \sum_{i=2}^q ({z^{(j)}})^TA_iz^{(j)} ({x^{(j)}})^TA_ix^{(j)}
+\rho^{(j)}\left(\|x^{(j)}-y^{(j)}\|^2+\|x^{(j)}-z^{(j)}\|^2 \right) \le \Upsilon,$
which by the new assumption on the boundedness of the level set $X_\Upsilon$ yields 
$
\rho^{(j)} \left(\|x^{(j)}-y^{(j)}\|^2+\|x^{(j)}-z^{(j)}\|^2 \right) \le \Upsilon - \min_{x\in X_\Upsilon} f(x),
$
which leads to $x^*=y^*=z^*$. Consequently, (i) and (ii) in Theorem \ref{thm: PPC convergence} hold for this case as well. 
\end{proof}

\end{document}